\documentclass[12pt]{article}
\usepackage{graphicx,setspace,amsmath,amssymb,subfigure,url,multirow,booktabs,verbatim,bm}
\usepackage{epsfig}
\usepackage{colortbl}
\usepackage[colorlinks,
            linkcolor=blue,
            anchorcolor=blue,
            citecolor=blue
            ]{hyperref}
\usepackage{listings}
\usepackage{mathrsfs}
\usepackage{amsmath,amssymb}
\usepackage{epstopdf}
\usepackage{algorithm}
\usepackage{algorithmic}
\usepackage{color}
\usepackage{algorithmic}
\usepackage[numbers,sort&compress]{natbib}

{\addtolength{\baselineskip}{.3\baselineskip}}
\def\ba{\begin{array}}
\def\ea{\end{array}}
\def\be{\begin{equation}}
\def\ee{\end{equation}}

\newtheorem{definition}{Definition}
\newtheorem{theorem}{Theorem}

\newtheorem{remark}{Remark}

\begin{document}

\begin{titlepage}
\title{{\large\bf Continuous grey model with conformable fractional derivative}}
\author{Wanli Xie$^{a\dag}$, Caixia Liu$^{a\dag}$,Weidong Li$^{a}$, Wenze Wu$^{b}$\thanks{corresponding author, E-mail: wenzew@mails.ccun.edu.cn; {\dag}: co-first authors.}, Chong Liu$^{c}$\\
  {\small $^{a}$Institute of EduInfo Science and Engineering, Nanjing Normal University, Nanjing Jiangsu 210097, China}\\
  {\small $^{b}$School of Economics and Business Administration, Central China Normal University, Wuhan 430079, China}\\
  {\small $^{c}$School of  Science,Inner Mongolia Agricultural University,Hohhot 010018,China}\\
    \\
}

\date{}

\maketitle\thispagestyle{empty}
\begin{abstract}
The existing fractional grey prediction models mainly use discrete fractional-order difference and accumulation, but in the actual modeling, continuous fractional-order calculus has been proved to have many excellent properties, such as hereditary. Now there are grey models established with continuous fractional-order calculus method, and they have achieved good results. However, the models are very complicated in the calculation and are not conducive to the actual application. In order to further simplify and improve the grey prediction models with continuous fractional-order derivative, we propose a simple and effective grey model based on conformable fractional derivatives in this paper, and two practical cases are used to demonstrate the validity of the proposed model.
\end{abstract}

\\

\end{titlepage}

\vspace{5mm}

\thispagestyle{empty}
\newpage

\pagestyle{plain}

\pagenumbering{arabic}

\section{Introduction}\label{sec:intro}

Fractional calculus has been around for hundreds of years and came around the same time as classical calculus. After years of development, fractional calculus has been widely applied in control theory, image processing, elastic mechanics, fractal theory, energy, medicine, and other fields \cite{1,2,3,4,5,6,7,8}. Fractional calculus is an extension of the integer-order calculus and the common fractional derivatives include Grunwald-Letnikov (GL) \cite{9}, Riemann-Liouville (RL) \cite{10}, Caputo \cite{11}, and so on. Although continuous fractional-order grey models have been applied in various fields, it is seldom used in the grey systems, while discrete fractional-order difference is mostly used at present.

The grey model was first proposed by Professor Deng. It solves the problem of small sample modeling, and the grey model does not need to know the distribution rules of data \cite{GM}. The potential rules of data can be fully mined through sequence accumulation, which has a broad application \cite{GM}. With the development of grey theory during several decades, grey prediction models have been developed very quickly and have been applied to all walks of life. For example, Li et al. \cite{13} used a grey prediction model to predict building settlements. Cao et al. \cite{14} proposed a dynamic Verhulst model for commodity price and demand prediction. Zhang et al. \cite{15} applied a grey prediction model and neural network model for stock prediction. Ma et al. \cite{16} presented a multi-variable grey prediction model for China's tourism revenue forecast. Wu et al. \cite{17} proposed a fractional grey Bernoulli model to forecast the renewable energy of China. Zeng et al. \cite{18} used a new grey model to forecast natural gas demand. Wu et al. \cite{19} put forward a fractional grey model for air quality prediction. Ding et al. \cite{20} presented a multivariable grey model for the prediction of carbon dioxide in China. Modeling background in the real world becomes more and more complex, which puts forward higher requirements for modeling. Many scholars have improved various grey prediction models. For example, Xie et al. \cite{21} proposed a grey model and the prediction formula was derived directly from the difference equation, which improved the prediction accuracy. Cui et al. \cite{22} presented a grey prediction model and it can fit an inhomogeneous sequence, which improved the range of application of the model. Chen et al. \cite{23} put forward a nonlinear Bernoulli model, which can capture nonlinear characteristics of data. Wu et al. \cite{FGM} proposed a fractional grey prediction model and it successfully extended the integer-order to the fractional-order, at the same time, they proved that the fractional-order grey model had smaller perturbation bounds integer order derivative. Ma et al. \cite{25} put forward a fractional-order grey prediction model that was simple in the calculation and was easy to be popularized and applied in engineering. Zeng et al. \cite{26} proposed an adaptive grey prediction model based on fractional-order accumulation. Wei et al. \cite{27} presented a method for optimizing the polynomial model and obtained expected results. Liu et al. \cite{28} proposed a grey Bernoulli model based on the Weibull Cumulative Distribution, which improved the modeling accuracy. In \cite{29}, a mathematical programming model was established to optimize the parameters of grey Bernoulli.

Although the above models have achieved good results, they all use continuous integer-order derivatives. In fact, the continuous derivative has many excellent characteristics, such as heritability \cite{30}. At present, there is little work on the grey prediction model based on continuous fractional derivative, and the corresponding research is still in early stage. In recent years, a new limit-based fractional order derivative is introduced by Khalil et al. in 2014 \cite{CF_define}, which is called ¡®¡®the conformable fractional derivative¡¯¡¯. It is simpler than the previous fractional order derivatives, such as the Caputo derivative and Riemann-Liouville derivative, so it can easily solve many problems, compared with other derivatives with complex definitions. In 2015, Abdeljawad \cite{develop_com} developed this new fractional order derivative and proposed many very useful and valuable results, such as Taylor power series expansions, Laplace transforms based on this novel fractional order derivative. Atangana et al. \cite{properties_cf} introduced the new properties of conformable derivative and proved some valuable  theorems. In 2017, Al-Rifae and Abdeljawad proposed \cite{Complexity_com} a regular fractional generalization of the Sturm-Liouville eigenvalue problems and got some important results. The Yavuz and Ya\c{s}k\i ran \cite{equation_com} suggested a new method for the approximate-analytical solution of the fractional one-dimensional cable differential equation (FCE) by employing the conformable fractional derivative. In this paper, we propose a new grey model based on conformable fractional derivative, which has the advantage of simplicity and efficiency.
The organization of this paper is as follows:

In the second section, we introduce a few kinds of fractional-order derivatives. In the third section, we show a grey model with Caputo fractional derivative and in the fourth section, we present a new grey prediction model containing conformable derivative. In the fifth section, we give the optimization methods of the order and background-value coefficient. In the sixth section, two practical cases are used to verify the validity of the model and the seventh section is a summary of the whole paper.

\section{Fractional-order derivative}
Fractional derivatives have rich forms, three common forms are Grunwald-Letnikov (GL), Riemann-Liouville (RL), and Caputo \cite{31}.

\begin{definition}[See \cite{31}]GL derivative with $\alpha$ order of function $f(t)$ is defined as
\begin{equation}
_a^{GL}D_t^\alpha f(t) = \sum\limits_{k = 0}^n {\frac{{{f^{(k)}}(a){{(t - a)}^{ - \alpha  + k}}}}{{\Gamma ( - \alpha  + k + 1)}}}  + \frac{1}{{\Gamma (n - \alpha  + 1)}}\int_a^t {{{(t - \tau )}^{n - \alpha }}} {f^{(n + 1)}}(\tau )d\tau
\end{equation}
where $_a^{GL}D_t^\alpha$ is the form of fractional derivative of GL, $\alpha  > 0, n - 1 < \alpha  < n, n \in N$, $[a,t]$ is the integral interval of $f(t)$, $\Gamma ( \cdot )$ is Gamma function, which has the following properties: $\Gamma (\alpha ) = \int_0^\infty  {{t^{\alpha  - 1}}} {e^{ - t}}dt$.
\end{definition}
\begin{definition}[See \cite{31}]
RL derivative with order $\alpha$ of function $f(t)$ is defined as

\begin{equation}
_a^{RL}D_t^\alpha f(t){\rm{ }} = {\frac{{{d^n}}}{{d{t^n}}}_a}D_t^{ - (n - \alpha )}f(t) = \frac{1}{{\Gamma (n - \alpha )}}\frac{{{d^n}}}{{d{t^n}}}\int_a^t {{{(t - \tau )}^{n - \alpha  - 1}}} f(\tau )d\tau
\end{equation}
where $_a^{RL}D_t^\alpha f(t)$ is the fractional derivative of RL, $a$ is an initial value, $\alpha$ is the order, $\Gamma (\cdot)$ is Gamma function.
\end{definition}
\begin{definition}[See \cite{31}]
Caputo derivative with $\alpha$-order of function $f(t)$ is defined as

\begin{equation}
_a^{C}D_t^\alpha f(t) = \frac{1}{{\Gamma (n - \alpha )}}\int_a^t {{{(t - \tau )}^{n - \alpha  - 1}}} {f^{(n)}}(\tau )d\tau ,{\rm{Among  them,}}_a^{C}D_t^\alpha f(t)
\end{equation}
where $a$ is an initial value, $\alpha$ is the order, $\Gamma (\cdot)$ is Gamma function. In particular, if the derivative order is ranged from 0 to 1, the Caputo derivative can be written as follows
\begin{equation}
_a^{C}D_t^\alpha f(t) = \frac{1}{{\Gamma (1 - \alpha )}}\int_a^t {{{(t - \tau )}^{ - \alpha }}} {f^\prime }(\tau )d\tau
\end{equation}
\end{definition}
Although the above derivatives have been successfully applied in various fields, it is difficult to be applied in engineering practice due to the complicated definition in the calculation. In recent years, some scholars have proposed a simpler fractal derivative called conformable derivative \cite{32} defined as follows.
\begin{definition}[See \cite{32}]
Assume ${T_\alpha }(f)(t)$ is the derivative operator of $f$: $[0,\infty) \to R$, $t > 0$, $\alpha  \in (0,1)$, and ${T_\alpha }(f)(t)$ is defined as

\begin{equation}
{T_\alpha }(f)(t) = \mathop {\lim }\limits_{\varepsilon  \to 0} \frac{{f\left( {t + \varepsilon {t^{1 - \alpha }}} \right) - f(t)}}{\varepsilon }
\end{equation}

when $\alpha  \in (n,n + 1]$, $f$ is differentiable at $t (t >0)$, the $\alpha$-order derivative of the function $f$ is

\begin{equation}
{T_\alpha }(f)(t) = \mathop {\lim }\limits_{\varepsilon  \to 0} \frac{{{f^{(\lceil \alpha  \rceil - 1)}}\left( {t + \varepsilon {t^{(\lceil \alpha  \rceil - \alpha )}}} \right) - {f^{(\lceil \alpha  \rceil - 1)}}(t)}}{\varepsilon }
\end{equation}
where $\lceil \alpha  \rceil$ is the smallest integer greater than or equal to $\alpha$.
\end{definition}
The conformable derivative has the following properties,
\begin{definition}[See \cite{32}]
Let $\alpha  \in (0,1]$ and $f$, $g$ be $\alpha$-differentiable at a point $t > 0$, then

(1) ${T_\alpha}(af + bg) = a{T_\alpha }(f) + b{T_\alpha }(g)$ for all $a, b \in R$.

(2) ${T_\alpha}\left( {{t^p}} \right) = p{t^{p - \alpha }}$ for all $p \in R$.

(3) ${T_\alpha}(\lambda) = 0$, for all constant functions $f(t) = \lambda$.

(4) ${T_\alpha} (fg) = f{T_\alpha }(g) + g{T_\alpha }(f)$.

(5) ${T_\alpha}\left( {\frac{f}{g}} \right) = \frac{{g{T_\alpha }(f) - f{I_\alpha }(g)}}{{{g^2}}}$.
where $T_\alpha$  is $a$-order conformable derivative.
\end{definition}
\begin{theorem}[See \cite{32}]
Let $\alpha  \in (0,1]$  and $f$, $g$  be $\alpha$-differentiable at a point $t > 0$. Then

\begin{equation}
\label{pro}
{T_\alpha }(f)(t) = {t^{1 - \alpha }}\frac{{df}}{{dt}}(t)
\end{equation}
\end{theorem}

{\it\textbf{ Proof.}} Let $h = \varepsilon {t^{1 - \alpha}}$, then $\begin{array}{l}
{T_\alpha }(f)(t) = \mathop {\lim }\limits_{\varepsilon  \to 0} \frac{{f\left( {t + \varepsilon {t^{1 - \alpha }}} \right) - f(t)}}{\varepsilon }= {t^{1 - \alpha }}\mathop {\lim }\limits_{h \to 0} \frac{{f(t + h) - f(t)}}{h}= {t^{1 - \alpha }}\frac{{df(t)}}{{dt}}
\end{array}$,

where $\frac{{df}}{{dt}}$is first-order Riemann derivative, ${T_\alpha }(f)(t)$ is $a$-order conformable derivative.

\begin{definition}[See \cite{32}] $I_\alpha ^a(f)(t) = I_1^a\left( {{t^{\alpha  - 1}}f} \right) = \int_a^t {\frac{{f(x)}}{{{x^{1 - \alpha}}}}} dx$, where the integral is the usual Riemann improper integral, and $\alpha  \in (0,1)$.

Based on the above definitions, we give the definitions of conformable fractional-order difference and derivative.
\end{definition}
\begin{definition}[See \cite{33}] The conformable fractional accumulation (CFA) of $f$ with $\alpha$-order is expressed as
\begin{equation}
\label{cfade}
\begin{array}{l}
{\nabla ^\alpha }f(k) = \nabla \left( {{k^{\alpha  - 1}}f(k)} \right) = \sum\limits_{i = 1}^k {\frac{{f(i)}}{{{i^{1 - \alpha }}}}} ,\alpha  \in (0,1],k \in {N^ + }\\
{\nabla ^\alpha }f(k) = {\nabla ^{(n + 1)}}\left( {{k^{\alpha  - [\alpha ]}}f(k)} \right),\alpha  \in (n,n + 1],k \in {N^ + }.
\end{array}
\end{equation}

The conformable fractional difference (CFD) of  $f$ with $\alpha$-order is given by
\begin{equation}
\begin{array}{*{20}{l}}
{{\Delta ^\alpha }f(k) = {k^{1 - \alpha }}\Delta f(k) = {k^{1 - \alpha }}[f(k) - f(k - 1)],\alpha  \in (0,1],k \in {N^ + }}\\
{{\Delta ^\alpha }f(k) = {k^{[\alpha ] - \alpha }}{\Delta ^{n + 1}}f(k),\alpha  \in (n,n + 1],k \in {N^ + }}.
\end{array}
\end{equation}
\end{definition}
In the next section, We give a brief introduce for the fractional grey model with Caputo derivative. This model uses continuous fractional derivative for modeling at the first time and achieves good results.
\section{Grey model with Caputo fractional derivative}
Most of the previous grey models were based on integer-order derivatives. Wu first proposed a grey prediction model based on the Caputo fractional derivative, and the time response sequence of the model was directly derived from the fractional derivative of Caputo, which achieved good results \cite{34}. In this section, we will introduce the modeling mechanism of this model.
\begin{definition}[See \cite{34}]  Assume ${X^{(0)}} = \left\{ {{x^{(0)}}(1),{x^{(0)}}(2), \cdots ,{x^{(0)}}(n)} \right\}$ is a non-negative sequence, the grey model with univariate of $p  (0 < p < 1)$ order equation is
\begin{equation}
\label{gmp_b1}
{\alpha ^{(1)}}{x^{(1 - p)}}(k) + a{z^{(0)}}(k) = b
\end{equation}
\end{definition}
where ${z^{(0)}}(k) = \frac{{{x^{(1 - p)}}(k) + {x^{(1 - p)}}(k - 1)}}{2}$, ${\alpha ^{(1)}}{x^{(1 - p)}}(k)$ is a $p$-order difference of ${x^{(0)}}(k)$, the least squares estimation of $GM(p,1)$ satisfies
$\left[ {\begin{array}{*{20}{l}}
a\\
b
\end{array}} \right] = {\left( {{B^{\rm{T}}}B} \right)^{ - 1}}{B^{\rm{T}}}Y$, where

\begin{equation}
B = \left[ {\begin{array}{*{20}{c}}
{ - {z^{(0)}}(2)}&1\\
{ - {z^{(0)}}(3)}&1\\
 \vdots & \vdots \\
{ - {z^{(0)}}(n)}&1
\end{array}} \right],Y = \left[ {\begin{array}{*{20}{c}}
{{\alpha ^{(1)}}{x^{(1 - p)}}(2)}\\
{{\alpha ^{(1)}}{x^{(1 - p)}}(3)}\\
 \vdots \\
{{\alpha ^{(1)}}{x^{(1 - p)}}(n)}
\end{array}} \right]
\end{equation}

The winterization equation of $GM(p,1)$ is
\begin{equation}
\label{GMP}
\frac{{{{\rm{d}}^p}{x^{(0)}}(t)}}{{{\rm{d}}{t^p}}} + a{x^{(0)}}(t) = b.
\end{equation}

Assume ${\hat x^{(0)}}(1) = {x^{(0)}}(1) $, the solution of the fractional equation calculated by the Laplace transform is

\begin{equation}
{x^{(0)}}(t) = \left( {{x^{(0)}}(1) - \frac{b}{a}} \right)\sum\limits_{k = 0}^\infty  {\frac{{{{\left( { - a{t^p}} \right)}^k}}}{{\Gamma (pk + 1)}}}  + \frac{b}{a}
\end{equation}

Then, the restored values of can be obtained

\begin{equation}
\label{resGMP}
{x^{(0)}}(k) = \left( {{x^{(0)}}(1) - \frac{b}{a}} \right)\sum\limits_{i = 0}^\infty  {\frac{{{{\left( { - a{k^p}} \right)}^i}}}{{\Gamma (pi + 1)}}}  + \frac{b}{a}
\end{equation}
Although many fractional grey models have achieved good results, most of the fractional gray prediction models use fractional difference and fractional accumulation, while those model still use integer derivative. Although there are some studies on grey models with fractional derivatives, they are more complicated to calculate than previous grey models. In order to simplify  calculation, we will propose a novel fractional prediction model with conformable derivative.
\section{Grey system model with conformable fractional derivative}

In this section, based on the conformable derivative, we propose a simpler grey model, named continuous conformable fractional grey model, abbreviated as CCFGM(1,1). Wu et al. \cite{35} first gives the unified form of conformable fractional accumulation Eq. (\ref{cfade}). On this basis, we use the matrix method to give the equivalent form of unified conformable fractional order accumulation.
\begin{theorem}The conformable fractional accumulation is
\label{CFA}
\begin{equation}
{x^{(\alpha )}}(k) = \sum\limits_{i = k}^n {\left[ {\begin{array}{*{20}{c}}
{ \lceil \alpha  \rceil }\\
{k - i}
\end{array}} \right]} \frac{{{x^{(0)}}(i)}}{{i \lceil \alpha  \rceil  - \alpha }},\alpha  \in {R^ + }
\end{equation}
where $\lceil \alpha  \rceil$  is the smallest integer greater than or equal to $\alpha$, $\left[ {\begin{array}{*{20}{c}}
{ \lceil \alpha  \rceil }\\
{k - i}
\end{array}} \right] = \frac{{ \lceil \alpha  \rceil ( \lceil \alpha  \rceil  + 1) \cdots ( \lceil \alpha  \rceil  + i - 1)}}{{\left( {k - i} \right)!}} = \left( {\begin{array}{*{20}{c}}
{k - i +  \lceil \alpha  \rceil  - 1}\\
{k - i}
\end{array}} \right) = \frac{{(k - i +  \lceil \alpha  \rceil  - 1)!}}{{(k - i)!( \lceil \alpha  \rceil  - 1)!}}$. $\alpha$ is the order of the model. Theoretically, the order of grey model can be any positive number. In order to simplify the calculation, we will make the order of the model between 0 and 1 in the later modeling.
\end{theorem}
{\it\textbf{ Proof.}} if $\alpha  \in (0,1], \lceil \alpha  \rceil  = 1$,

$\begin{array}{l}
{x^{(\alpha )}}(k) = \sum\limits_{i = 1}^k {\frac{{{x^{(0)}}(i)}}{{{i^{1 - \alpha }}}}} {x^{(0)}}(i) = \left[ {{x^{(0)}}(1),{x^{(0)}}(2), \cdots ,{x^{(0)}}(n)} \right]\left( {\begin{array}{*{20}{c}}
1&1& \cdots &1&1\\
0&{\frac{1}{{{2^{1 - \alpha }}}}}& \cdots &{\frac{1}{{{2^{1 - \alpha }}}}}&{\frac{1}{{{2^{1 - \alpha }}}}}\\
 \vdots & \vdots & \vdots & \vdots & \vdots \\
0&0& \cdots &{\frac{1}{{{{\left( {n - 1} \right)}^{1 - \alpha }}}}}&{\frac{1}{{{{\left( {n - 1} \right)}^{1 - \alpha }}}}}\\
0&0& \cdots &0&{\frac{1}{{{n^{1 - \alpha }}}}}
\end{array}} \right)\\
{\rm{          }}\\
 = \left[ {{x^{(0)}}(1),{x^{(0)}}(2), \cdots ,{x^{(0)}}(n)} \right]\left( {\begin{array}{*{20}{c}}
1&0& \cdots &0&0\\
0&{\frac{1}{{{2^{1 - \alpha }}}}}& \cdots &0&0\\
 \vdots & \vdots & \vdots & \vdots & \vdots \\
0&0& \cdots &{\frac{1}{{{{\left( {n - 1} \right)}^{1 - \alpha }}}}}&0\\
0&0& \cdots &0&{\frac{1}{{{n^{1 - \alpha }}}}}
\end{array}} \right)\left( {\begin{array}{*{20}{c}}
1&1& \cdots &1&1\\
0&1& \cdots &1&1\\
 \vdots & \vdots & \vdots & \vdots & \vdots \\
0&0& \cdots &1&1\\
0&0& \cdots &0&1
\end{array}} \right)\\
 = \sum\limits_{i = 1}^k {\left( {\begin{array}{*{20}{c}}
{k - i}\\
{k - i}
\end{array}} \right)} \frac{{x(i)}}{{{i^{1 - \alpha }}}}$, $k = 1,2, \cdots, n.
\end{array}$

if $\alpha  \in (1,2], \lceil \alpha  \rceil = 2$.

$\begin{array}{l}
{x^{(\alpha)}}(j) = \sum\limits_{j = k}^n {\sum\limits_{i = j}^n {\frac{{{x^{(0)}}(i)}}{{{i^{2 - r}}}}} } \\
 = {\left[ \begin{array}{l}
{x^{(0)}}(1)\\
{x^{(0)}}(2)\\
 \cdots \\
{x^{(0)}}(n)
\end{array} \right]^{\rm{T}}}\left( {\begin{array}{*{20}{c}}
1&0& \cdots &0&0\\
0&{\frac{1}{{{2^{2 - \alpha }}}}}& \cdots &0&0\\
 \vdots & \vdots & \vdots & \vdots & \vdots \\
0&0& \cdots &{\frac{1}{{{{\left( {n - 1} \right)}^{2 - \alpha }}}}}&0\\
0&0& \cdots &0&{\frac{1}{{{n^{2 - \alpha }}}}}
\end{array}} \right)\left( {\begin{array}{*{20}{c}}
1&1& \cdots &1&1\\
0&1& \cdots &1&1\\
 \vdots & \vdots & \vdots & \vdots & \vdots \\
0&0& \cdots &1&1\\
0&0& \cdots &0&1
\end{array}} \right)\left( {\begin{array}{*{20}{c}}
1&1& \cdots &1&1\\
0&1& \cdots &1&1\\
 \vdots & \vdots & \vdots & \vdots & \vdots \\
0&0& \cdots &1&1\\
0&0& \cdots &0&1
\end{array}} \right)\\
 = {\left[ \begin{array}{l}
{x^{(0)}}(1)\\
{x^{(0)}}(2)\\
 \cdots \\
{x^{(0)}}(n)
\end{array} \right]^{\rm{T}}}\left( {\begin{array}{*{20}{c}}
1&0& \cdots &0&0\\
0&{\frac{1}{{{2^{2 - \alpha }}}}}& \cdots &0&0\\
 \vdots & \vdots & \vdots & \vdots & \vdots \\
0&0& \cdots &{\frac{1}{{{{\left( {n - 1} \right)}^{2 - \alpha }}}}}&0\\
0&0& \cdots &0&{\frac{1}{{{n^{2 - \alpha }}}}}
\end{array}} \right)\left( {\begin{array}{*{20}{c}}
1&2& \cdots &{n - 1}&n\\
0&1& \cdots &{n - 2}&{n - 1}\\
 \vdots & \vdots & \vdots & \vdots & \vdots \\
0&0& \cdots &1&2\\
0&0& \cdots &0&1
\end{array}} \right)\\
 = {\left[ \begin{array}{l}
{x^{(0)}}(1)\\
{x^{(0)}}(2)\\
 \cdots \\
{x^{(0)}}(n)
\end{array} \right]^{\rm{T}}}\left( {\begin{array}{*{20}{c}}
1&0& \cdots &0&0\\
0&{\frac{1}{{{2^{2 - \alpha }}}}}& \cdots &0&0\\
 \vdots & \vdots & \vdots & \vdots & \vdots \\
0&0& \cdots &{\frac{1}{{{{\left( {n - 1} \right)}^{2 - \alpha }}}}}&0\\
0&0& \cdots &0&{\frac{1}{{{n^{2 - \alpha }}}}}
\end{array}} \right)\left( {\begin{array}{*{20}{c}}
1&{\left( {\begin{array}{*{20}{c}}
2\\
1
\end{array}} \right)}& \cdots &{\left( {\begin{array}{*{20}{c}}
{n - 1}\\
{n - 2}
\end{array}} \right)}&{\left( {\begin{array}{*{20}{c}}
n\\
{n - 1}
\end{array}} \right)}\\
0&1& \cdots &{\left( {\begin{array}{*{20}{c}}
{n - 2}\\
{n - 3}
\end{array}} \right)}&{\left( {\begin{array}{*{20}{c}}
{n - 1}\\
{n - 2}
\end{array}} \right)}\\
 \vdots & \vdots & \vdots & \vdots & \vdots \\
0&0& \cdots &1&{\left( {\begin{array}{*{20}{c}}
2\\
1
\end{array}} \right)}\\
0&0& \cdots &0&1
\end{array}} \right)\\
 = \sum\limits_{i = 1}^k {\left( {\begin{array}{*{20}{c}}
{k - i + 1}\\
{k - i}
\end{array}} \right)} \frac{{x(i)}}{{{i^{2 - \alpha }}}}$, $k = 1,2, \cdots, n.
\end{array}$

Assuming that the equation holds when $\alpha  \in (m - 1,m]$, then $\lceil \alpha  \rceil =m$,
${x^{(\alpha )}}(k) = \sum\limits_{i = k}^n {\left[ {\begin{array}{*{20}{c}}
m\\
{k - i}
\end{array}} \right]} \frac{{{x^{(0)}}(i)}}{{i \lceil \alpha  \rceil  - \alpha }},\alpha  \in {R^ + }$, when $\alpha  \in (m,m+1]$, $\lceil \alpha  \rceil = m+1$,\\
let $\left( {\begin{array}{*{20}{c}}
1&0& \cdots &0&0\\
0&{\frac{1}{{{2^{m + 1 - \alpha }}}}}& \cdots &0&0\\
 \vdots & \vdots & \vdots & \vdots & \vdots \\
0&0& \cdots &{\frac{1}{{{{\left( {n - 1} \right)}^{m + 1 - \alpha }}}}}&0\\
0&0& \cdots &0&{\frac{1}{{{n^{m + 1 - \alpha }}}}}
\end{array}} \right){\rm{ = A}}$, we have

$\begin{array}{l}
{x^\alpha }(k) = {\left[ \begin{array}{l}
{x^{(0)}}(1)\\
{x^{(0)}}(2)\\
 \cdots \\
{x^{(0)}}(n)
\end{array} \right]^{\rm{T}}}{\rm{A}}{\left[ {\begin{array}{*{20}{c}}
1&0& \cdots &0&0\\
1&1& \cdots &0&0\\
 \vdots & \vdots & \ddots & \vdots & \vdots \\
1&1& \cdots &1&0\\
1&1& \cdots &1&1
\end{array}} \right]^{m + 1}}\\
{\rm{  = }}{\left[ \begin{array}{l}
{x^{(0)}}(1)\\
{x^{(0)}}(2)\\
 \cdots \\
{x^{(0)}}(n)
\end{array} \right]^{\rm{T}}}{\rm{A}}\left( {\begin{array}{*{20}{c}}
1&{\left( {\begin{array}{*{20}{c}}
m\\
1
\end{array}} \right)}& \cdots &{\left( {\begin{array}{*{20}{c}}
{m + n - 2}\\
{n - 1}
\end{array}} \right)}\\
0&1& \cdots &{\left( {\begin{array}{*{20}{c}}
{m + n - 3}\\
{n - 2}
\end{array}} \right)}\\
 \vdots & \vdots & \vdots & \vdots \\
0&0& \cdots &{\left( {\begin{array}{*{20}{c}}
m\\
1
\end{array}} \right)}\\
0&0& \cdots &1
\end{array}} \right)\left( {\begin{array}{*{20}{c}}
1&1& \cdots &1&1\\
0&1& \cdots &1&1\\
 \vdots & \vdots & \vdots & \vdots & \vdots \\
0&0& \cdots &1&1\\
0&0& \cdots &0&1
\end{array}} \right)\\
 = {\left[ \begin{array}{l}
{x^{(0)}}(1)\\
{x^{(0)}}(2)\\
 \cdots \\
{x^{(0)}}(n)
\end{array} \right]^{\rm{T}}}{\rm{A}}\left( {\begin{array}{*{20}{c}}
1&{\left( {\begin{array}{*{20}{c}}
m\\
0
\end{array}} \right) + \left( {\begin{array}{*{20}{c}}
m\\
1
\end{array}} \right)}& \cdots &{\sum\limits_{i = 0}^{n - 3} {\left( {\begin{array}{*{20}{c}}
{m + i}\\
{i + 1}
\end{array}} \right)} }&{\sum\limits_{i = 0}^{n - 2} {\left( {\begin{array}{*{20}{c}}
{m + i}\\
{i + 1}
\end{array}} \right)} }\\
0&1& \cdots &{\sum\limits_{i = 0}^{n - 4} {\left( {\begin{array}{*{20}{c}}
{m + i}\\
{i + 1}
\end{array}} \right)} }&{\sum\limits_{i = 0}^{n - 3} {\left( {\begin{array}{*{20}{c}}
{m + i}\\
{i + 1}
\end{array}} \right)} }\\
 \vdots & \vdots & \vdots & \vdots & \vdots \\
0&0& \cdots &1&{\left( {\begin{array}{*{20}{c}}
m\\
0
\end{array}} \right) + \left( {\begin{array}{*{20}{c}}
m\\
1
\end{array}} \right)}\\
0&0& \cdots &0&{}
\end{array}} \right)\\
 = {\left[ \begin{array}{l}
{x^{(0)}}(1)\\
{x^{(0)}}(2)\\
 \cdots \\
{x^{(0)}}(n)
\end{array} \right]^{\rm{T}}}{\rm{A}}\left( {\begin{array}{*{20}{c}}
1&{\left( {\begin{array}{*{20}{c}}
{m + 1}\\
1
\end{array}} \right)}& \ldots &{\left( {\begin{array}{*{20}{c}}
{m + n - 2}\\
{n - 2}
\end{array}} \right)}&{\left( {\begin{array}{*{20}{c}}
{m + n - 1}\\
{n - 1}
\end{array}} \right)}\\
0&1& \cdots &{\left( {\begin{array}{*{20}{c}}
{m + n - 3}\\
{n - 3}
\end{array}} \right)}&{\left( {\begin{array}{*{20}{c}}
{m + n - 2}\\
{n - 2}
\end{array}} \right)}\\
 \vdots & \vdots & \vdots & \vdots & \vdots \\
0&0& \cdots &1&{\left( {\begin{array}{*{20}{c}}
{m + 1}\\
1
\end{array}} \right)}\\
0&0& \cdots &0&1
\end{array}} \right)\\
 = \sum\limits_{i = 1}^k {\left( {\begin{array}{*{20}{c}}
{k - i + m}\\
{k - i}
\end{array}} \right)} {x^{(0)}}(i)
\end{array}$.
So the result is proved.
\begin{remark}
 Similarly, the Refs. \cite{25,35} give the other two methods to get the same result. It can be proved that the definitions of these accumulation are essentially the same. Using the matrix method can help us better understand the fractional accumulation. Secondly, it can better help us write computer programs.
\end{remark}
Next, we will derive the grey differential equation with continuous conformable fractional derivatives.
\begin{definition}
Assume ${X^{(0)}} = \left\{ {{x^{(0)}}(1),{x^{(0)}}(2), \cdots ,{x^{(0)}}(n)} \right\}$ is a non-negative sequence, $r(0 < r < 1)$-order winterization equation can be dedined as follows,
\begin{equation}\label{CCFGM}
\frac{{{{\rm{d}}^r}{x^{(q)}}(t)}}{{{\rm{d}}{t^r}}} + a{x^{(q)}}(t) = b,
\end{equation}
where ${X^{(q)}} = \left( {{x^{(q)}}(1),{x^{(q)}}(2), \cdots ,{x^{(q)}}(n)} \right)$ is the q-order $(0<q<1)$ cumulative sequence of ${X^{(0)}}$, and $\frac{{{{\rm{d}}^r}{x^{(q)}}(t)}}{{{\rm{d}}{t^r}}} = {T_r}\left( {{x^{(q)}}(1)} \right)$ is continuous conformable fractional-order derivative.
\end{definition}

\begin{remark}
If r=1 and q=1, the equation (\ref{CCFGM}) is equivalent to GM(1,1) (see \cite{GM}), if $r \in \left[ {0,1} \right]$ and q =0, the equation  equation (\ref{CCFGM}) is equivalent to the  equation (\ref{GMP}) in form, if r=0 and $q \in \left[ {0,1} \right]$, the equation  equation (\ref{CCFGM}) is equivalent to the  FGM(1,1) (see \cite{FGM}) in form.
\end{remark}

\begin{theorem}
The exact solution of the conformable fractional-order differential equation is
\begin{equation}
\label{CCF_time_response}
{\hat x^{(r)}}(k) = \frac{{\hat b + \left( {\widehat a{x^{(0)}}(1) - \hat b} \right)e{}^{\frac{{\widehat a\left( {1 - {k^r}} \right)}}{r}}}}{{\widehat a}},k = 1,2,3,...,n{\rm{ }}(n > 4)
\end{equation}
\end{theorem}

{\it\textbf{Proof.}} Using equation (\ref{pro}) to convert the fractional order derivative into integer order derivative, we can find the exact solution of equation (\ref{CCFGM}).
$\frac{{{{\rm{d}}^r}{x^{(q)}}(t)}}{{{\rm{d}}{t^r}}} + a{x^{(q)}}(t) = b$,

${t^{1 - r}}\frac{{{\rm{d}}{x^{(q)}}(t)}}{{{\rm{d}}t}} + a{x^{(q)}}(t) = b$, by integrating the two sides, we have $\int {\frac{{{\rm{d}}{x^{(q)}}(t)}}{{(b{\rm{ - }}a{x^{(q)}}(t))}}}  = \int {\frac{{dt}}{{{t^{1 - r}}}}}$, so

$\ln \left| {b - a{x^{(q)}}(t)} \right| = \left( {\frac{{ - a}}{r}} \right){t^r} + {C_{\rm{1}}}$, $b - a{x^{(q)}}(t =  \pm {e^{{C_{\rm{1}}}}}e{}^{\left( {\frac{{ - a}}{r}} \right){t^r}}$, it can be sorted out,

${x^{(q)}}(t) = \frac{{b{\rm{ + }}Ce{}^{\left( {\frac{{ - a}}{r}} \right){t^r}}}}{a}$, assume $\hat a$, $\hat b$ is estimated parameters, ${\hat x^{(0)}}(k)$ is an estimated value of ${x^{(0)}}(k)$, $k$ is a discrete variable with respect to $t$, with ${\hat x^{(q)}}(0) = {x^{(0)}}(1)$, then $C = \left( {\widehat a{x^{(0)}}(1) - \hat b} \right){e^{\frac{{\widehat a}}{r}}}$, so the time response function of the CCFGM model is Eq. (\ref{CCF_time_response}).

\begin{remark}
If r=1 and q=1, the equation (\ref{CCF_time_response}) is equivalent to response function of GM(1,1) (see \cite{GM}), if $r \in \left[ {0,1} \right]$ and q =0, the equation  equation (\ref{CCF_time_response}) is equivalent to the  equation (\ref{resGMP}) in form (Mittag Leffler is a direct generalization of exponential function.), if r=0 and $q \in \left[ {0,1} \right]$, the equation  equation (\ref{CCF_time_response}) is equivalent to the response function of FGM(1,1) (see \cite{FGM}) in form.
\end{remark}
  Next, we will derive the discrete form of CCFGM(1,1) model. Through the discrete difference equation, we can use least squares algorithm to get the parameters of the model.
The predicted value can be obtained by q-order difference of the obtained predicted value, as follows,
${{\hat x}^{(0)}}(k) = \Delta {\nabla ^{1 - q}}{{\hat x}^{(q)}}(k)$.
${x^{(1 - q)}}(t)$ stands for $1-q$-order accumulation, and it is equal to $\Delta {}^1{\nabla ^q}{x^{(0)}}(t)$. ${\nabla ^q}{x^{(0)}}(t)$ is the q-order accumulation of ${x^{(0)}}(t)$, $\Delta {}^r{x^{(q)}}(t)$ is the r-order difference of ${x^{(q)}}(t)$, $q \in \left[ {0,1} \right]$.
\begin{theorem}
The difference equation of the continuous conformable grey model is
\begin{equation}
\label{dccfgm}
{x^{(q - r)}}(t) + a\frac{{\rm{1}}}{{\rm{2}}}\left[ {{x^{(q)}}(k - 1) + {x^{(q)}}(k)} \right] = b, q \in \left[ {0,1} \right],  r \in \left[ {0,1} \right].
\end{equation}
\end{theorem}
{\it\textbf{ Proof.}} Integrate CCFGM with $r$-order on both sides of Eq. (\ref{CCFGM}):
\begin{equation}
\iint  \cdots \int_{k - 1}^k {\frac{{{d^r}{x^{(q)}}}}{{d{t^r}}}} d{t^r} + a\iint  \cdots \int_{k - 1}^k {{x^{(q)}}} (t)d{t^r} = b\iint  \cdots \int_{k - 1}^k d {t^r}
\end{equation}
where
\begin{equation}\label{dccfgm1}
\iint  \cdots \int_{k - 1}^k {\frac{{{d^r}{x^{(q)}}(t)}}{{d{t^r}}}} d{t^r} \approx \frac {\Delta ^r}{x^{(q)}} = {x^{(q - r)}}(t)
\end{equation}
${x^{(q - r)}}(t)$ stands for $q-r$-order accumulation, and it is equal to $\Delta {}^r{\nabla ^q}{x^{(0)}}(t)$. ${\nabla ^q}{x^{(0)}}(t)$ is the q-order accumulation of ${x^{(0)}}(t)$, $\Delta {}^r{x^{(0)}}(t)$ is the r-order difference of ${x^{(0)}}(t)$, $r \in \left[ {0,1} \right]$,
$r \in \left[ {0,1} \right]$.

According to the generalized trapezoid formula (see \cite{mao}), we have,
\begin{equation}\label{dccfgm2}
\iint  \cdots \int_{k - 1}^k {{x^{(q)}}} (t)d{t^r} \approx \frac{{\rm{1}}}{{\rm{2}}}\left[ {{x^{(q)}}(k - 1) + {x^{(q)}}(k)} \right]
\end{equation}
According to equation (\ref{gmp_b1}) and equation (\ref{GMP}),we have
\begin{equation}\label{dccfgm3}
\iint  \cdots \int_{k - 1}^k b d {t^r}=b \iint  \cdots \int_{k - 1}^k d {t^r}\approx \int_{k - 1}^k {bd} t \approx b.
\end{equation}
By equation (\ref{dccfgm1}), equation (\ref{dccfgm2}), and equation (\ref{dccfgm3}), the basic form of CCFGM(1,1) can be written as equation (\ref{dccfgm}).

Through the least square method, we can get the parameter of the CCFGM(1,1) is
\begin{equation}
\label{parameter_model}
\hat{a} = \left[ {\begin{array}{*{20}{l}}
a\\
b
\end{array}} \right] = {\left( {{B^{\rm{T}}}B} \right)^{ - 1}}{B^{\rm{T}}}Y
\end{equation}
where
\begin{equation}\label{B_para}
B = \left[ {\begin{array}{*{20}{c}}
{ - \frac{{\rm{1}}}{{\rm{2}}}\left[ {{x^{(q)}}(1) + {x^{(q)}}(2)} \right]}&1\\
{ - \frac{{\rm{1}}}{{\rm{2}}}\left[ {{x^{(q)}}(2) + {x^{(q)}}(3)} \right]}&1\\
 \vdots & \vdots \\
{ - \frac{{\rm{1}}}{{\rm{2}}}\left[ {{x^{(q)}}(n - 1) + {x^{(q)}}(n)} \right]}&1
\end{array}} \right],Y = \left[ {\begin{array}{*{20}{c}}
{{x^{_{(q - r)}}}(2)}\\
{{x^{_{(q - r)}}}(3)}\\
 \vdots \\
{{x^{_{(q - r)}}}(n)}
\end{array}} \right]
\end{equation}

Let  $\varepsilon  = Y - B\hat a$  be the error sequence and $s=\varepsilon \cdot \varepsilon^{\mathrm{T}}=Y-B \hat{a}^{\mathrm{T}}(Y-B \hat{a})=\sum\limits_{k = 2}^n {{{\left\{ {{x^{(q - r)}}(t) + a\frac{1}{2}\left[ {{x^{(q)}}(k - 1) + {x^{(q)}}(k)} \right] - b} \right\}}^2}} dx$,
when s is minimized, values of a and b satisfy
\begin{equation}
\left\{ {\begin{array}{*{20}{l}}
{\frac{{\partial s}}{{\partial a}} = \sum\limits_{k = 2}^n {\left\{ {{x^{(q - r)}}(t) + a\frac{1}{2}\left[ {{x^{(q)}}(k - 1) + {x^{(q)}}(k)} \right] - b} \right\}} \left[ {{x^{(q)}}(k - 1) + {x^{(q)}}(k)} \right]dx = 0}\\
{\frac{{\partial s}}{{\partial b}} =  - 2\sum\limits_{k = 2}^n {\left\{ {{x^{(q - r)}}(t) + a\frac{1}{2}\left[ {{x^{(q)}}(k - 1) + {x^{(q)}}(k)} \right] - b} \right\}}  = 0}
\end{array}} \right.,
\end{equation}
where ${\hat a}$ is defined in the Eq. (\ref{parameter_model}), $B$ and $Y$ defined in the Eq. (\ref{B_para}).
\section{Optimization of the optimal difference order $r$ and accumulation order $q$}
The accumulative order is usually given by default, but in fact, the difference order $r$ and accumulation order $q$ as part of the model greatly affect the model accuracy. Their values can be dynamically adjusted according to different modeling content. So the correct order of the model are particularly important. In the following, we first established the following mathematical programming model to optimize the two super parameters and used a whale optimization algorithm for optimization \cite{36}.

\begin{equation}
\begin{array}{l}
{\min _{r,q}}\frac{1}{n}\sum\limits_{i = 1}^n {\left| {\frac{{{{\hat x}^{(0)}}\left( {{k_i}} \right) - {x^{(0)}}\left( {{k_i}} \right)}}{{{x^{(0)}}\left( {{k_i}} \right)}}} \right|}  \times 100\% \\
{\rm{s}}{\rm{.t}}\left\{ \begin{array}{l}
r \in \left[ {0,1} \right],q \in \left[ {0,1} \right]\\
{x^{(q)}}(k) = \sum\limits_{i = k}^n {\left[ {\begin{array}{*{20}{c}}
{ \lceil q \rceil }\\
{k - i}
\end{array}} \right]{x^{(0)}}(i)} \frac{{x(i)}}{{{i^{ \lceil q \rceil  - q}}}},q > 0\\
B = \left[ {\begin{array}{*{20}{c}}
{ - \frac{{\rm{1}}}{{\rm{2}}}\left[ {{x^{(q)}}(1) + {x^{(q)}}(2)} \right]}&1\\
{ - \frac{{\rm{1}}}{{\rm{2}}}\left[ {{x^{(q)}}(2) + {x^{(q)}}(3)} \right]}&1\\
 \vdots & \vdots \\
{ - \frac{{\rm{1}}}{{\rm{2}}}\left[ {{x^{(q)}}(n - 1) + {x^{(q)}}(n)} \right]}&1
\end{array}} \right],Y = \left[ {\begin{array}{*{20}{c}}
{{x^{_{(q - r)}}}(2)}\\
{{x^{_{(q - r)}}}(3)}\\
 \vdots \\
{{x^{_{(q - r)}}}(n)}
\end{array}} \right]\\
{{\hat x}^{(q)}}(k) = \frac{{\hat b + \left( {\widehat a{x^{(0)}}(1) - \hat b} \right)e{}^{\frac{{\widehat a\left( {1 - {k^r}} \right)}}{r}}}}{{\widehat a}},k = 2,3,4,...,n{\rm{ }}(n > 4)\\
{{\hat x}^{(0)}}(k) = \Delta {\nabla ^{1 - q}}{{\hat x}^{(q)}}(k)
\end{array} \right.{\rm{ }}
\end{array}
\end{equation}

\section{Application}

In order to verify the validity of the model, we test the model with two actual cases, and compare it with other forecasting models.

\textbf{Case 1.} Prediction of domestic energy consumption in China (Ten thousand ton standard coal)

{In this case, we select the data of domestic energy consumption in China from 2005 to 2015 for fitting and the data from 2016 to 2017 for testing. The corresponding results are shown in Table \ref{tcase1} and Figure 1.}

\begin{figure}[H]
\begin{minipage}[t]{0.5\linewidth}
\centering
\includegraphics[height=5cm,width=8cm]{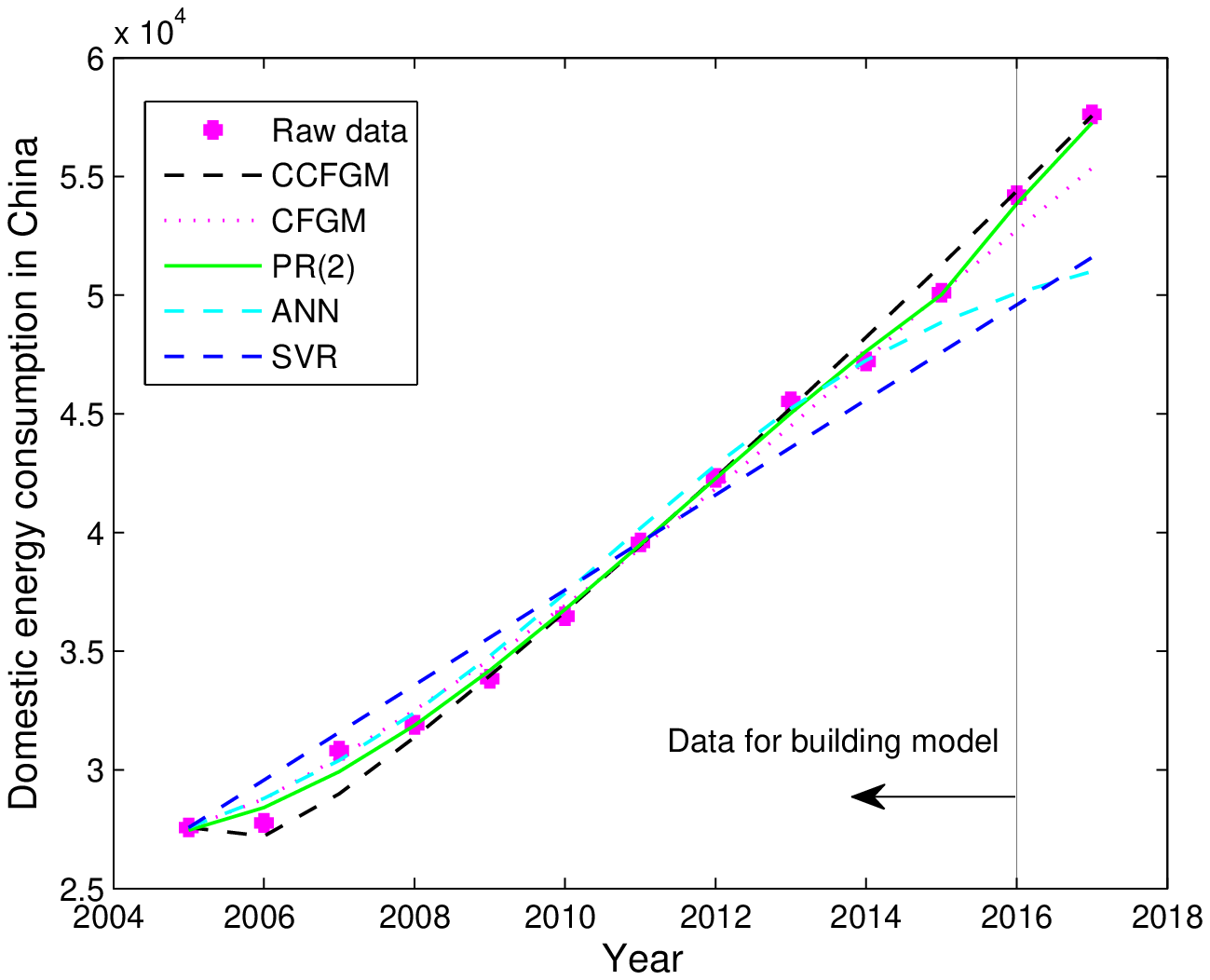}
\caption{Test results of five models.}
\end{minipage}
\hfill
\begin{minipage}[t]{0.5\linewidth}
\centering
\includegraphics[height=5cm,width=8cm]{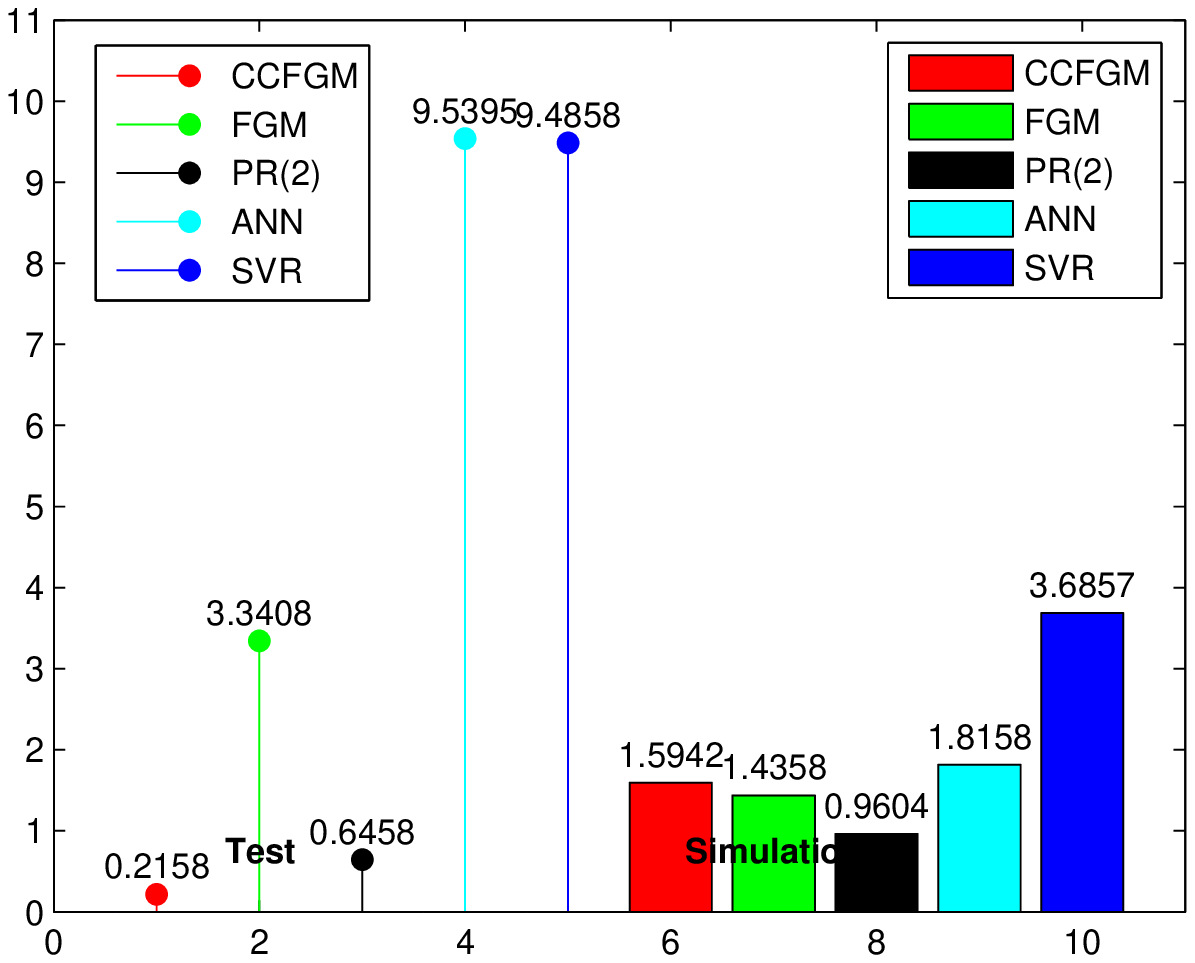}
\caption{Error comparison of five grey models.}
\end{minipage}
\label{fc1}
\end{figure}

\begin{table}[H]
\caption{Comparison of test results of five grey models.}\centering
\begin{tabular}{p{12mm}p{18mm}p{18mm}p{18mm}p{18mm}p{18mm}p{25mm}}
\hline
Year	&Raw data	&FGM	&PR(2)	     &ANN	    &SVR	    &CCFGM\\
\hline
2005	&27573	&27573.00 	&27461.01 	&27576.24 	&27572.90 	&27573.00\\
2006	&27765	&28776.69 	&28414.85 	&28801.47 	&29574.73 	&27207.68\\
2007	&30814	&30529.07 	&29920.42 	&30403.47 	&31576.57 	&28992.48\\
2008	&31898	&32510.82 	&31879.33 	&32409.59 	&33578.40 	&31373.76\\
2009	&33843	&34650.97 	&34193.17 	&34790.27 	&35580.23 	&33965.07\\
2010	&36470	&36925.23 	&36763.53 	&37441.64 	&37582.07 	&36671.07\\
2011	&39584	&39324.40 	&39492.00 	&40193.56 	&39583.90 	&39459.34\\
2012	&42306	&41845.55 	&42280.18 	&42848.65 	&41585.73 	&42317.23\\
2013	&45531	&44488.93 	&45029.65 	&45235.79 	&43587.57   &45239.62\\
2014	&47212	&47256.57 	&47642.01 	&47249.65 	&45589.40 	&48224.63\\
2015	&50099	&50151.64 	&50018.85 	&48859.33 	&47591.23 	&51271.90\\
\hline
MAPE &          & 1.4358 	&0.9604 	&1.8158     &3.6857 	&1.5942\\
\hline
2016	&54209	&52721.73 	&53852.78 	&50091.35 	&49593.07 	&54381.79\\
2017	&57620	&55350.93 	&57254.38 	&51003.38 	&51594.90 	&57555.00\\
\hline
MAPE &		   & 3.3408 	&0.6458 	 &   9.5395 &9.4858 	&0.2158\\
\hline
\end{tabular}
\label{tcase1}
\end{table}
{The test errors of five grey models are shown in Figure 2. The experimental results show that the fitting error and test error of the proposed model are 1.5942\% and 0.2158\% respectively, and the fitting error and test error of the FGM model are 1.4358\% and 3.3408\% respectively. The fitting error and test error of PR(2) are 0.9604\% and 0.6458\%, respectively, ANN are 1.8158\% and 9.5395\%, respectively, SVR are 3.6857\% and 9.4858\%, respectively. The fitting errors of PR(2) are slight lower than ours. However, the test error of our model are smaller than other models.}

\textbf{Case 2.} Prediction of domestic coal consumption in China (ten thousand tons).
{Coal consumption is related to the sustainable development of society. Accurate and effective prediction of coal consumption can contribute to effective decision-making and early warning.}

\begin{figure}[H]
\begin{minipage}[t]{0.5\linewidth}
\centering
\includegraphics[height=5cm,width=8cm]{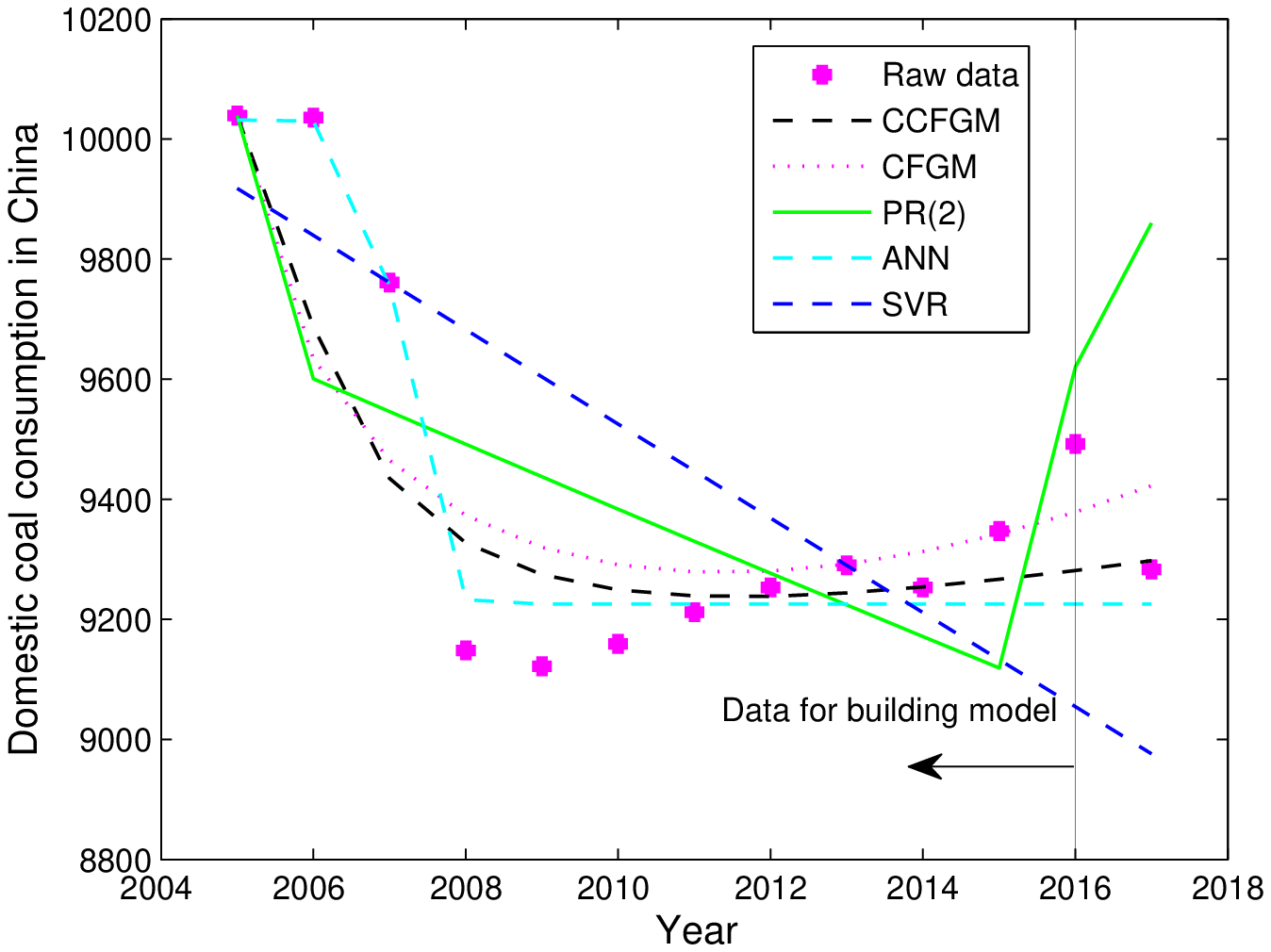}
\caption{Test results of five models.}
\end{minipage}
\hfill
\begin{minipage}[t]{0.5\linewidth}
\centering
\includegraphics[height=5cm,width=8cm]{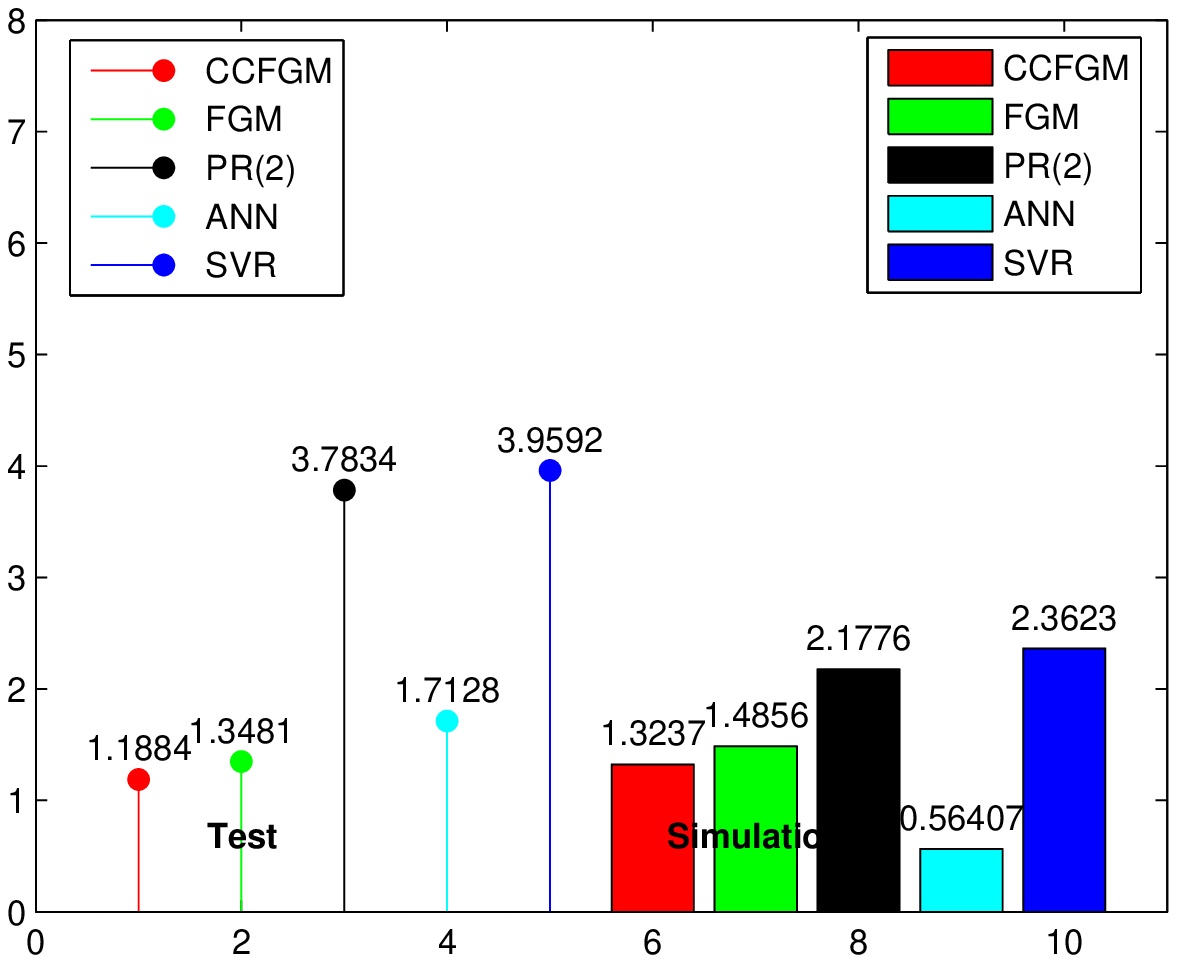}
\caption{Error comparison of five grey models.}
\end{minipage}
\end{figure}

\begin{table}[H]
\caption{Comparison of test results of five grey models.}\centering
\begin{tabular}{p{12mm}p{18mm}p{18mm}p{18mm}p{18mm}p{18mm}p{25mm}}
\hline
Year	&Raw data	&FGM	&PR(2)	     &ANN	    &SVR	    &CCFGM\\
\hline
2005	&10039.00 	&10039.00 	&10039.00 	&10031.67 	&9917.90 	&10039.00\\
2006	&10036.00 	&9633.62 	&9600.63 	&10029.69 	&9839.40 	&9687.54\\
2007	&9761.00 	&9464.86 	&9545.90 	&9755.68 	&9760.90 	&9434.76\\
2008	&9148.00 	&9372.94 	&9491.48 	&9232.88 	&9682.40 	&9326.97\\
2009	&9122.00 	&9319.89 	&9437.36 	&9225.60 	&9603.90 	&9274.31\\
2010	&9159.00 	&9290.99 	&9383.56 	&9225.57 	&9525.40 	&9248.88\\
2011	&9212.00 	&9279.07 	&9330.07 	&9225.57 	&9446.90 	&9238.87\\
2012	&9253.00 	&9280.18 	&9276.87 	&9225.57 	&9368.40 	&9238.35\\
2013	&9290.00 	&9291.95 	&9223.99 	&9225.57 	&9289.90 	&9243.99\\
2014	&9253.00 	&9312.86 	&9171.40 	&9225.57 	&9211.40 	&9253.82\\
2015	&9347.00 	&9341.96 	&9119.11 	&9225.57 	&9132.90 	&9266.55\\
\hline
MAPE &              &1.4856 	&2.1776 	&0.5641 	&2.3623 	&1.3237\\
2016	&9492.00 	&9378.60 	&9620.06 	&9225.57 	&9054.40 	&9281.34\\
2017	&9283.00 	&9422.38 	&9860.18 	&9225.57 	&8975.90 	&9297.62\\
\hline
MAPE &              &1.3481 	    &3.7834 	&1.7128 	&3.9592 	&1.1884\\
\hline
\end{tabular}
\label{tcase2}
\end{table}
{Table \ref{tcase2}, Figure 3 and Figure 4 show the prediction of carbon dioxide emission with our model. From Table \ref{tcase2}, we can see that the fitting error and test error of our model are 1.3237\% and 1.1884\%, respectively. The fitting error and test error of FGM model are 1.4856\% and 1.3481\%, respectively, PR(2) are 2.1776\% and 3.7834\%, respectively, ANN are 0.5641\% and 1.7128\%, respectively, SVR are 2.3623\% and 3.9592\%, respectively. It can be seen that  our model has smaller test errors compared with other models, which means that our model is superior to other models.
}

\section{Conclusion}

In this paper, we propose a grey forecasting model with a conformable fractional derivative. Compared with integer derivatives, continuous fractional derivatives have been proved to have many excellent properties. However, the most existing grey models are modeled by integer derivatives. Secondly, it has been proved that the integer derivative cannot simulate some special development laws in nature, the model can be further optimized by extending the grey model with the integer derivative to the fractional derivative. The existing fractional order grey model with continuous fractional-order derivative, achieved good result, but its calculation is complicated. This paper proposes a new grey model with conformable fractional-order derivative, further to simplify the calculation. Two actual cases show that our model has high precision, and it can be easily promoted in engineering. The contributions of this paper are as follows:

(1) We constructed a fractional-order differential equation with a conformable derivative as a whitening form of our model.

(2) We built the mathematical programming model to optimize the order and of CCFGM(1,1) by whale optimizer, which further improved the prediction accuracy of the model.

(3) We verify the validity of the model in this paper through two actual cases. This model with a simpler structure can achieve similar or even better accuracy than other models.

Although the model in this paper has some advantages, it can be further improved from the following aspects:

(1) In order to improve the modeling accuracy of the model, a more efficient optimization algorithm can be used to optimize parameters.

(2) The model proposed in this paper is linear and cannot capture the nonlinear characteristics of the data. Accordingly, nonlinear characteristics can be studied for establishing a more universal and robust grey prediction model.

\section {Conflicts of Interest}
No potential conflict of interest was reported by the authors.

\section {Acknowledgement}

The work was supported by grants from the Postgraduate Research \& Practice Innovation Program of Jiangsu Province [Grant No.KYCX19\_0733]; grants from the Postgraduate Research \& Practice Innovation Program of Jiangsu Province [Grant No.KYCX20\_1144].


\begin{thebibliography}{}


\bibitem{1}
Abdeljawad T, Al-Mdallal QM. Discrete Mittag-Leffler kernel type fractional difference initial value problems and Gronwalls inequality. J Comput Appl Math 2018;339:218-30.

\bibitem{2}
Abdeljawad T, Al-Mdallal QM, Jarad F. Fractional logistic models in the frame of fractional operators generated by conformable derivatives. Chaos Solit Fract 2019;119:94-101.

\bibitem{3}
Al-Mdallal QM. On fractional-Legendre spectral Galerkin method for fractional Sturm-Liouville problems. Chaos Solit Fract 2018;116:261-7.

\bibitem{4}
Al-Mdallal QM, Hajji MA. A convergent algorithm for solving higher-order nonlinear fractional boundary value problems. Fract Calcul Appl Anal 2015;18(6):1423-40.

\bibitem{5}
Al-Mdallal Q, Abro KA, Khan I. Analytical solutions of fractional Walters b fluid with applications. Complexity 2018:1-11.

\bibitem{6}
Al-Mdallal QM, Hajji MA. A convergent algorithm for solving higher-order nonlinear fractional boundary value problems. Fract Calcul Appl Anal 2015;18(6):1423-40.

\bibitem{7}
Al-Mdallal QM, Omer ASA. Fractional-order legendre-collocation method for solving fractional initial value problems. Appl Math Comput 2018;321:74-84.

\bibitem{8}
Al-Mdallal QM, Syam MI. An efficient method for solving non-linear singularly perturbed two points boundary-value problems of fractional order. Commun Nonlinear Sci Numer Simul 2012;17(6):2299-308.

\bibitem{9}
Scherer R, Kalla S L, Tang Y, et al. The Gr¡§1nwald¡§CLetnikov method for fractional differential equations. Computers \& Mathematics with Applications 2011;62(3):902-917.

\bibitem{10}
WANG, Bing. S-asymptotically -periodic Solutions of R-L Fractional Derivative-Integral Equation. Technology Vision  2015;(17):155-155.

\bibitem{11}
Wang J, Lv L, Zhou Y. Ulam stability and data dependence for fractional differential equations with Caputo derivative. Electronic Journal of Qualitative Theory of Differential Equations 2011;63(63):1-10.

\bibitem{GM}
Liu S, Yang Y, Forrest J. Grey Data Analysis: Methods, Models and Applications, Springer Press, 2016.


\bibitem{13}
LI R, Wang L, ZHANG S. Study of the application of grey forecast model in settlement forecast of high buildings. Journal of Earth Science and Enivronmental 2005;27(1):84-87.

\bibitem{14}
Cao C, Gu X. The application of grey dynamic Verhulst metabolism model in the prediction of products' prices and demands. INFORMATION AND CONTROL-SHENYANG- 2005;34(4):398.

\bibitem{15}
Zhang Q, Zhu H. Application of grey model and neural network in stock prediction. computer engineering \& applications  2013;49(12):242-236.

\bibitem{16}
Xin M, Zhibin L, Yong W. Application of a novel nonlinear multivariate grey Bernoulli model to predict the tourist income of China. Journal of Computational \& Applied Mathematics 2018;347:84-94.

\bibitem{17}
Wu W, Ma X, Zeng B, et al. Forecasting short-term renewable energy consumption of China using a novel fractional nonlinear grey Bernoulli model. Renewable energy 2019;140:70-87.


\bibitem{18}
Zeng B, Li C. Prediction the natural gas demand in China using a self-adapting intelligent grey model. Energy 2016;112:810-825.

\bibitem{19}
Wu L, Zhao H. Forecasting Air Quality Indicators for 33 Cities in China. CLEAN¡§CSoil, Air, Water  2020;48(1):1-12.

\bibitem{20}
Ding S, Dang Y G, Li X M, et al. Prediction Chinese CO2 emissions from fuel combustion using a novel grey multivariable model.Journal of Cleaner Production  2017;162:1527-1538.

\bibitem{21}
Xie N M, Liu S F. Discrete grey prediction model and its optimization. Applied Mathematical Modelling 2009;33(2):1173-1186.

\bibitem{22}
Cui J, Liu S F, Zeng B, et al. A novel grey prediction model and its optimization. Applied Mathematical Modelling 2013;37(6):4399-4406.

\bibitem{23}
Chen C I. Application of the novel nonlinear grey Bernoulli model for forecasting unemployment rate. Chaos, Solitons \& Fractals  2008;37(1):278-287.

\bibitem{FGM}
Wu L, Liu S, Yao L, et al. Grey system model with the fractional order accumulation. Communications in Nonlinear Science \& Numerical Simulation 2013;18(7):1775-1785.

\bibitem{25}
Ma X, Wu W, Zeng B, et al. The conformable fractional grey system model. ISA transactions  2020;96:255-271.

\bibitem{26}
Zeng B, Liu S. A self$-$ adaptive intelligence gray prediction model with the optimal fractional order accumulating operator and its application. Mathematical Methods in the Applied Sciences  2017;40(18):7843-7857.

\bibitem{27}
Wei B, Xie N, Hu A. Optimal solution for novel grey polynomial prediction model. Applied Mathematical Modelling  2018;62:717-727.

\bibitem{28}
Liu X, Xie N. A nonlinear grey forecasting model with double shape parameters and its application. Applied Mathematics and Computation  2019;360:203-212.

\bibitem{29}
Wang Z X, Hipel K W, Wang Q, et al. An optimized NGBM(1.1) model for prediction the qualified discharge rate of industrial wastewater in China. Applied Mathematical Modelling 2011;35(12):5524-5532.

\bibitem{30}
Zhao D, Pan X, Luo M. A new framework for multivariate general conformable fractional calculus and potential applications. Physica A: Statistical Mechanics and its Applications  2018;510:271-280.

\bibitem{CF_define}
Khalil Roshdi, Al Horani M, Yousef Abdelrahman, Sababheh Mohammad. A new definition of fractional derivative. J Comput Appl Math
2014;264:65¨C70.

\bibitem{develop_com}
Abdeljawad T. On conformable fractional calculus. Journal of computational and Applied Mathematics 2015;279:57-66.

\bibitem{properties_cf}
Atangana A, Baleanu D, Alsaedi A. New properties of conformable derivative. Open Mathematics  2015;13(1):889-898.

\bibitem{Complexity_com}
Al-Refai M, Abdeljawad T. Fundamental results of conformable Sturm-Liouville eigenvalue problems Complexity;2017:1-8.

\bibitem{equation_com}
Yavuz M, Ya\c{s}k\i ran B.Approximate-analytical solutions of cable equation using conformable fractional operator. New Trends in Mathematical Sciences 2017;5(4):209-219.

\bibitem{31}
Wang J, Wen Y, Gou Y, et al. Fractional-order gradient descent learning of BP neural networks with Caputo derivative. Neural Networks  2017;89:19-30.

\bibitem{32}
Khalil R, Horani M A, Yousef A, et al. A new definition of fractional derivative. Journal of Computational \& Applied Mathematics  2014;264(5):65-70.

\bibitem{33}
Ma X, Wu W, Zeng B, et al. The conformable fractional grey system model. ISA transactions  2020;96:255-271.

\bibitem{34}
Wu L F, Liu S F, Yao L G. Grey model with Caputo fractional order derivative. System Engineering$-$Theory \& Practice  2015;35(5):1311-1316.

\bibitem{35}
Wu W, Ma X, Zhang Y, et al. A novel conformable fractional non-homogeneous grey model for forecasting carbon dioxide emissions of BRICS countries. Science of the Total Environment  2020;707:135447.
\bibitem{mao}
Mao S , Gao M , Xiao X , et al. A novel fractional grey system model and its application[J]. Applied Mathematical Modelling, 2016, 40(7-8):5063-5076.

\bibitem{36}
Gene H, Charles C. Matrix computations. The Johns Hopkins University Press, Baltimore and London, 1996.


\end{thebibliography}
\end{document}